\documentclass[12pt]{article}
  
 \def\lddots{\mathinner{\mkern1mu\raise1pt\hbox{.}\mkern2mu  
  
\raise4pt\hbox{.}\mkern2mu\raise7pt\vbox{\kern7pt\hbox{.}}\mkern1mu}}  
\makeatletter  
\def\numberbysection{\@addtoreset{equation}{section}  
 \def\theequation{\thesection.\arabic{equation}}}  
  
\makeatother  
  
\numberbysection

\newcommand{\be}{\begin{eqnarray}}  
\newcommand{\ee}{\end{eqnarray}}  
\newcommand{\non}{\nonumber}  
  
\newcommand{\tr}{\mathop{\rm tr}\nolimits}  
  
\begin{document}  
  
  
\strut\hfill{LAPTH-974/03}  
\begin{center}  
  
\LARGE Commuting quantum traces: the case of reflection algebras \\[0.8in]  
\large  { Jean Avan} \footnote{e-mail: avan@ptm.u-cergy.fr}\\  
\normalsize{{Laboratory of Theoretical Physics and Modelization\\ University of Cergy,  
5 mail Gay--Lussac, Neuville--sur--Oise,\\F-95031 Cergy--Pontoise Cedex}}\\  
\vspace{.2in}  
\large  {Anastasia Doikou} \footnote{e-mail: doikou@lapp.in2p3.fr} \\  
\normalsize{Theoretical Physics Laboratory of Annecy--Le--Vieux,\\  
LAPTH, B.P. 110, Annecy--Le--Vieux, F-74941, France}\\  
  
\end{center}

\begin{abstract}  
We formulate a systematic construction of commuting quantum  
traces for reflection algebras. This is achieved by introducing two  
dual sets of generalized reflection equations with associated  
consistent fusion procedures. Products of their respective solutions yield commuting 
quantum traces.  
  
\end{abstract}  
\section{Introduction}  
  
The concept of quantum traces which will be discussed here  goes  
back to the work of Maillet \cite{MAI1} where commuting quantum  
analogues of the classical Poisson-commuting traces of powers of  
Lax matrices $TrL^{n}, \Big \{ n \in \textbf{N} \Big \}$ were  
explicitly constructed in the context of quantum group structures.  
Starting from the well known fundamental quantum group relation  
\be R_{12}(\lambda_{1} - \lambda_{2})\ L_{1q}(\lambda_{1})\  
L_{2q}(\lambda_{2}) = L_{2q}(\lambda_{2})\ L_{1q}(\lambda_{1})\ R_{1  
2}(\lambda_{1} - \lambda_{2})\,, \label{YBEf} \ee where the  
quantum $R$ matrix obeys the Yang Baxter equation, \cite{baxter,  
korepin} \be R_{12}(\lambda_{1} - \lambda_{2})\ R_{1  
3}(\lambda_{1})\ R_{23}(\lambda_{2}) = R_{23}(\lambda_{2})\ R_{1  
3}(\lambda_{1})\ R_{1 2}(\lambda_{1} - \lambda_{2})\,, \label{YBE}  
\ee (and $q$ denotes the quantum space \footnote{the quantum space may have the structure of a tensor product of single quantum spaces.} on which the generators of 
the quantum group inside the Lax matrix, act), quantum commuting  
objects were built  with the generic form \be H_{N} = Tr_{1 \ldots  
N}({\cal R}_{1\ldots N} L_{1q} L_{2q} \ldots L_{Nq} ). \ee ${\cal  
R}$ is any auxiliary-space operator which satisfies a commuting  
form of the Yang--Baxter equation: precisely  it is required that  
${\cal R}_{12 \ldots  N}$  commute with coproduct--like structures  
of the form $R_{1a}R_{2a}\ldots R_{Na}$. In particular, the choice  
${\cal R}_{12 \ldots N} = P_{12} R_{12} P_{23} R_{23}\ldots P_{N-1  
N}R_{N-1N}$, with $P_{ab}$ the permutation operator acting on the  
auxiliary spaces $V_{a} \otimes V_{b}$, leads to the exact quantum  
analogue of the classical traces $Tr( L^{N})$.  A similar possibility occurs
here justifying the characterization of our operators
as ``quantum traces''. This concept is different from the one leading 
to monodromy matrices, whose classical limit instead yields $(Tr L)^{N}$ and
powers thereof.
We shall give a brief account of the construction in \cite{MAI1}  
as a preliminary illustration of the mechanisms involved in our  
more complicated procedure.  
  
We approach here the problem of formulating  
a quantum trace construction for the more general quadratic algebras,  
discussed e.g. in \cite{MAI2}, of the form  
\be A_{12}\ T_{1}\ B_{12}\ T_{2} =  
T_{1}\ C_{12}\ T_{2}\ D_{12}, \label{braid} \ee where $T$ encapsulates generators  
of the algebra and $A$, $B$, $C$, $D$ are c-number structure matrices.
  
We shall focus here on a  specific set of quadratic algebras
known as ``reflection  
algebras'' \cite{cherednik} where the matrices $A$, $B$, $C$, $D$ are  
related to a single $R$ matrix, albeit a priori depending on one complex  
spectral parameter. This structure is defined as \cite{cherednik, sklyanin}:  
\be R_{12}(\lambda_{1} -\lambda_{2})\ T_{1} (\lambda_{1})\ R_{21}(\lambda_{1} +\lambda_{2})\ T_{2}(\lambda_{2})=T_{2}(\lambda_{2})\  R_{12}(\lambda_{1} +\lambda_{2})\ T_{1} (\lambda_{1})\ R_{21}(\lambda_{1} -\lambda_{2}).\ \label{re} \ee  
$ T_{i}(\lambda_{i})$ is understood as a matrix acting  
on a finite-dimensional auxiliary space $V_i$ with matrix entries being  
operators, representing the quantum reflection algebra, acting on some  
Hilbert space of quantum states hereafter denoted by the single label $q$.
  
Treatment of the most general case will be left for a forthcoming investigation
\cite{ADN} but it must be emphasized here that the particular case of reflection
algebras actually embraces the full algebraic richness of the general case.  
Explicit realizations of quantum $T$ matrices may easily be built
using the fundamental results in \cite{sklyanin}
from quantum Lax matrices $L_{q}(\lambda)$  obeying (\ref{YBE}) and  
quantum reflection matrices $K(\lambda)$ obeying (\ref{re})
 such as were constructed in \cite{AACFR} (see also references therein) as  
\be  
T_q(\lambda) = L_{q}(\lambda)K(\lambda)(L_{q}(\lambda))^{-1}.  
\ee  Note that the $K$ matrix may act on some extra boundary quantum space $q_{b}$ along the lines described in \cite{bako, dm}.
  
Furthermore it is known \cite{sklyanin,MAI2} that  
\be  
t_{qq'}(\lambda) =Tr_1 \Big  (K_{1q'}^{+}(\lambda)T_{1q}(\lambda)\Big )  
\label{trace1}\ee realize a family of mutually commuting operators on the quantum  
 Hilbert space $q \otimes q'$:  
\be \Big [t_{qq'}(\lambda),\ t_{qq'}(\mu) \Big ] =0,\label{com0}\ee where the
 matrix $K_{q'}^{+}$  
is a solution of the so called ``dual'' reflection equation \cite{MAI2,MENE}:  
\be && R_{12}(-\lambda_{1} +\lambda_{2})\ K_{1}^{+} (\lambda_{1})\ M_{1}^{-1}\ 
R_{21}(-\lambda_{1} -\lambda_{2}-2\rho)\ M_{1}\ K_{2}^{+}(\lambda_{2})\non\\ &=&  
K_{2}^{+}(\lambda_{2})\  M_{1}\ R_{12}(-\lambda_{1} -\lambda_{2}-2\rho)\  
M_{1}^{-1}\ K_{1}^{+} (\lambda_{1})\  
R_{21}(-\lambda_{1} +\lambda_{2})\,. \label{red} \ee Again 
$K_{i}^{+}$ is understood as a matrix acting  
on the same finite-dimensional auxiliary space $V_i$ with matrix entries
being generically operators, representing the quantum reflection algebra, acting on
a different Hilbert space of quantum states denoted by the single label $q'$. Furthermore the $K^{+}$ matrix may act on some boundary quantum space denoted by the index $q'_{b}$.

This notion  
of duality will be one crucial ingredient of our construction.  
  
The problem is to obtain extensions of (\ref{trace1}) as traces  
involving products (as opposed to {\it tensor} products) of more than one $T$ (monodromy) matrix. Results on quantum traces obtained  
for algebra (\ref{re}) will soon be extended \cite{ADN} to 
the generic quadratic  
Yang--Baxter algebras by a suitable generalization of the procedure
described here.  
  
Throughout this paper we shall impose  
several conditions on the $R$ matrix. We assume that:  
  
{\bf I.} It obeys the Yang--Baxter equation (\ref{YBE}).  
 
{\bf II.} It obeys the following symmetry requirement,  
\be  
R_{12}(\lambda) = R_{21}(\lambda)^{t_{1}t_{2}}. \label{transp}\ee For instance Yangian matrices constructed in \cite{AACFR} obey this condition.  
 
{\bf III.} It obeys unitarity and crossing relations \be  
R_{12}(\lambda)\ R_{21}(-\lambda)\propto 1, ~~R_{12}(\lambda)= V_{1}\ R_{12}^{t_{2}}(-\lambda-\rho)\ V_{1}=  
V_{2}^{t_{2}}\ R_{12}^{t_{1}}(-\lambda-\rho)\ V_{2}^{t_{2}}, \label{cross} \ee 
using (\ref{transp}) to obtain the second crossing relation.
 

Furthermore the $R$ matrix obeys the crossing--unitarity relation, \be  
R_{21}(\lambda)^{t_{1}}\ M_{1}^{-1}\ R_{12}(-\lambda  
-2\rho)^{t_{1}}\ M_{1} \propto 1, \ee 
and the commutation relation
\be 
\Big [R_{12},\ M_{1} M_{2} \Big ] =0.
\label{crosscom1}
\ee
 Here $V$ is a c-number  
matrix such that $V^{2} =1$ and $V^{t}\ V =M$. The crossing 
used in \cite{AACDFR} is of such type, and the Yangian $R$  
matrices constructed in \cite{AACFR} obey these conditions.

These conditions are particular formulations,  
in the context of reflection algebras,  
of more generic requirements for the structural  
matrices of  quadratic Yang--Baxter algebras (\ref{braid}).Let us immediately indicate that  general consistency 3-space exchange
conditions generalizing (\ref{YBE}) have already  been formulated,
see e.g. equations $(15)$ in \cite{MAI2}.  
  
Our results are formulated as two basic steps:  
  
{\bf 1.} We introduce two sets of dual  generalized  
reflection equations extending (\ref{re}), (\ref{red}),  
closely related to the fusion procedure described in e.g. \cite{MENE}. We  
establish the existence and the form of the simplest solutions in two basic Lemmas. More general solutions may then be constructed
 using a dressing procedure.  
  
{\bf 2.} We then build quantum traces for reflection algebras by combining  
solutions of the two dual sets of reflection equations in a form similar to  
(\ref{trace1}) with suitable non-trivial dressing operators.

We finally conclude with some comments on the perspectives  
opened by our results.  
  
\section{Commuting traces}  
  
\subsection{Review on Yang--Baxter type algebras}  
  
Before we describe the construction of commuting traces associated to reflection  
algebras, we first would like to present the construction of  
quantum traces associated to Yang--Baxter algebras, realized in \cite{MAI2}, from a slightly different point of view. Namely, we shall introduce a set of generalized  fundamental equations with the help of which we will be able to build the commuting  
 traces, and which we will use later in reflection algebras.  
  
Let us first fix some convenient notations based on the coproduct structure  
of equation (\ref{YBEf}) \cite{Dr,KRS}. We introduce the objects, hereafter denoted as  
``fused'' $R$ matrices,  
\be  
R_{{\cal N} {\cal M}'}(\lambda_{{\cal N}} - \lambda_{{\cal M}'}) &=&  
R_{11'}(\lambda_{1} -\lambda_{1'}) R_{21'}(\lambda_{2} -\lambda_{1'}) \ldots  
R_{n1'}(\lambda_{n} -\lambda_{1'}) \non\\  
&&  R_{12'}(\lambda_{1} -\lambda_{2'}) R_{22'}(\lambda_{2} -\lambda_{2'})  
\ldots  R_{n2'}(\lambda_{n} -\lambda_{2'}) \ldots \non\\  
&&  R_{1m'}(\lambda_{1} -\lambda_{m'}) R_{2m'}(\lambda_{2}  
-\lambda_{m'}) \ldots  R_{nm'}(\lambda_{n} -\lambda_{1m'}) \label{genr} \ee  
where one defines {\rm ordered} sets  
${\cal N}\equiv <12 \ldots n>$, ${\cal M}' \equiv <1'2'  
\ldots m'>$. We should emphasize that the notation we use here is  
essentially inspired by the one introduced in \cite{MENE}  
describing the fusion procedure for open spin chains. The above object  
satisfies the following properties: \be R_{{\cal N} {\cal M}'}(\lambda_{{\cal N}}- \lambda_{{\cal M}'})^{t_{{\cal M}'}t_{{\cal N}}} = R_{\bar {\cal M}' \bar {\cal N}}(\lambda_{\bar {\cal N}}- \lambda_{\bar {\cal M}'}), \label{p1} \ee  
and the generalization of crossing (\ref{cross}), 
\be R_{\bar {\cal N} {\cal M}'}(\lambda)= V_{{\cal N}}\ R_{{\cal N} {\cal M}'}(-\lambda-\rho)^{t_{{\cal M}'}}\ V_{{\cal N}}, ~~ R_{{\cal N}\bar  {\cal M}'}(\lambda)= V_{{\cal M}'}^{t_{{\cal M}'}}\ R_{{\cal N} {\cal M}'}(-\lambda-\rho)^{t_{{\cal N}}}\ V_{{\cal M}'}^{t_{{\cal M}'}}. \ee  
It also satisfies unitarity and crossing unitarity, i.e.  
\be &&  
R_{{\cal N} {\cal M}'}(\lambda_{{\cal N}}- \lambda_{{\cal M}'})\  
R_{\bar {\cal M}'\bar {\cal N}}(\lambda_{\bar {\cal M}'}-  
\lambda_{\bar {\cal N}}) \propto 1, \non\\ && R_{{\cal N} {\cal  
M}'}(\lambda_{{\cal N}}+ \lambda_{{\cal M}'})^{t_{{\cal M}'}}\ M_{{\cal M}'}^{-1}\  
R_{\bar {\cal M}'\bar {\cal N}}(-\lambda_{\bar {\cal M}'}-  
\lambda_{\bar {\cal N}}-2\rho)^{t_{{\cal M}'}}\ M_{{\cal M}'} \propto 1,  
\label{cross2} \ee 
and generalized commutation relation
\be 
\Big [R_{{\cal N} {\cal M}'},\ M_{{\cal N}} M_{{\cal M}'} \Big ] =0.
\label{crosscom2}
\ee 
Here one defines  {\rm anti--ordered} sets $\bar {\cal M}'  
\equiv <m'(m'-1) \ldots  
1'>$, $\bar {\cal N} \equiv <n(n-1) \ldots 1>$.  
Now we introduce the set of generalized fundamental equations:  
\be  
R_{{\cal N} \bar {\cal M}'}(\lambda_{{\cal N}} - \lambda_{\bar  
{\cal M}'})\ {\cal L}_{{\cal N}}(\lambda_{{\cal N}})\ {\cal L}_{{\cal  
M}'}(\lambda_{{\cal M}'})= {\cal L}_{{\cal M}'}(\lambda_{{\cal  
M}'})\ {\cal L}_{{\cal N}}(\lambda_{{\cal N}})\ R_{{\cal N} \bar {\cal  
M}'}(\lambda_{{\cal N}} - \lambda_{\bar {\cal M}'}). \label{gybe}  
\ee All solutions ${\cal L}_{{\cal N}}$ of the above equation are  
actually good candidates for the construction of quantum commuting traces. In  
particular, the objects  
\be H_{{\cal N}} = Tr_{{\cal N}}{\cal  
L}_{{\cal N}}(\lambda_{{\cal N}}) \label{trace} \ee  
realize a family of commuting operators $\Big [H_{{\cal N}},\ H_{{\cal  
M}'} \Big ] =0$. An obvious  
solution of equation (\ref{gybe}) is  
\be {\cal L}_{{\cal  
N}}(\lambda_{{\cal N}}) = L_{1}(\lambda_{1}) \ldots  
L_{n}(\lambda_{n}) \label{sol} \ee  
 where $L_{i}$'s are Lax matrices obeying the fundamental equation (\ref{YBEf}).  
 Obviously this solution leads to trivially decoupled traces, therefore  
one defines  ``dressed'' solutions of  
the generalized fundamental equation (see \cite{MAI1, MAI2}),  
for instance:  
\be {\cal L}_{{\cal N}}(\lambda_{{\cal N}}) = \check  
R_{12}(\lambda_{1} -\lambda_{2}) \check R_{23}(\lambda_{2}  
-\lambda_{3})\ldots \check R_{n-1 n}(\lambda_{n-1} -\lambda_{n})  
L_{1}(\lambda_{1}) \ldots L_{n}(\lambda_{n}), \ee where $\check R_{12} \equiv  
P_{12}R_{12}$ and $P_{12} $ is the operator exchanging auxiliary spaces $1$ and $2$  
and spectral parameters $\lambda_1$ and $\lambda_2$ (this last  
property is not easy to actually realize and we shall comment later  
on the practical aspects of this realization). The objects $\check R$  
are characterized by their  
commutation relations with fused $R$-matrices:  
\be \check R_{12}(\lambda_{1} -\lambda_{2})\ R_{13}(\lambda_{1} -\lambda_{3})\  
R_{23}(\lambda_{2} -\lambda_{3}) =R_{13}(\lambda_{1} -\lambda_{3})\  
R_{23}(\lambda_{2} -\lambda_{3})\ \check R_{12}(\lambda_{1} -\lambda_{2}).  
\label{dress0} \ee The role of $\check R$'s is in this sense purely technical: they  dress the solutions (\ref{sol}) so that the traces  
(\ref{trace}) have non--trivial (non--decoupled) structure, but due to (\ref{dress0})  
they do not modify the exchange relations, which guarantee commutation of (\ref{trace}).  
We emphasize that {\cal any} object obeying the commutation  
relation (\ref{dress0}) is a good dressing operator;  $\check R$ is simply  
an easily constructed example of it.  
  
\subsection{Commuting traces related to reflection algebras}  
  
We now come to our main concern, which is  
the explicit construction of the commuting traces related to  
reflection algebras. For this purpose we shall introduce  
the notions of generalized reflection  
equations and duals thereof.  
  
We define, in analogy to the case related to Yang--Baxter  
algebras, the set of generalized reflection equations associated  
to the fused $R$-matrices (\ref{genr}) as: \be && {\cal T}_{{\cal  
N}q}(\lambda_{{\cal N}})\ R_{{\cal M}'{\cal N}}(\lambda_{ {\cal  
N}}+\lambda_{ {\cal M}'})\ {\cal T}_{{\cal M}'q}(\lambda_{{\cal  
M}'})\ R_{\bar {\cal N} {\cal M}'}(-\lambda_{\bar  {\cal  
N}}+\lambda_{{\cal M}'}) \non\\  & =& R_{ {\cal M}'\bar {\cal  
N}}(-\lambda_{\bar  {\cal N}}+\lambda_{ {\cal M}'})\ {\cal  
T}_{{\cal M}'q}(\lambda_{{\cal M}'})\ R_{{\cal N}  {\cal  
M}'}(\lambda_{ {\cal N}}+\lambda_{ {\cal M}'})\ {\cal T}_{{\cal  
N}q}(\lambda_{{\cal N}}), \label{greq} \ee where the objects 
${\cal T}_{{\cal N}q}$ are matrices acting on  
tensor products of auxiliary spaces indexed by ordered sets ${\cal  
N}$ as $V_{1} \otimes V_{2} \ldots \otimes V_{n}$ with operator  
entries acting on the ``bulk'' quantum space labelled by  
index $q$. In general they may also act on some extra boundary  
quantum space (see e.g. \cite{bako, dm}), which is denoted  
by the  
 index $q_{b}$.  
  
We also introduce the set of  generalized ``dual'' reflection  
equations, which has the following structure \be && {\cal  
K}_{{\cal N}q'}(\lambda_{{\cal N}})\ M_{{\cal M}'}\ R_{\bar {\cal  
M}'\bar{\cal N}}(-\lambda_{\bar {\cal N}}-\lambda_{\bar {\cal  
M}'}-2\rho)\ M_{{\cal M}'}^{-1}\ {\cal K}_{{\cal  
M}'q'}(\lambda_{{\cal M}'})\ R_{{\cal N} \bar {\cal M}'}(\lambda_{  
{\cal N}}-\lambda_{\bar {\cal M}'}) \non\\  & =& R_{\bar {\cal  
M}'{\cal N}}(\lambda_{ {\cal N}}-\lambda_{\bar {\cal M}'})\ {\cal  
K}_{{\cal M}'q'}(\lambda_{{\cal M}'})\ M_{{\cal M}'}^{-1}\  
R_{\bar{\cal N} \bar {\cal M}'}(-\lambda_{\bar {\cal  
N}}-\lambda_{\bar {\cal M}'}-2\rho)\ M_{{\cal M}'}\ {\cal  
K}_{{\cal N}q'}(\lambda_{{\cal N}}). \label{greqd} \ee ${\cal  
K}_{{\cal N}q'}$ are similarly  matrices acting on tensor products of  
auxiliary spaces indexed by ordered sets ${\cal N}$ as $V_{1}  
\otimes V_{2} \ldots \otimes V_{n}$ with operator entries acting  
on a a priori different  ``bulk'' quantum space labelled by  
index $q'$. In general they may act on some extra boundary  
quantum space as well (see e.g. \cite{bako, dm}), which is denoted  
by the index $q'_{b}$. We would like at this stage to point out the structural  
similarity between the generalized reflection equation (\ref{greq}), and the  
corresponding reflection equation for fused $K$ matrices introduced in \cite{MENE}.  In addition we assume the existence of a transposition antimorphism
$t_{q'q'_{b}}$ acting on the operator entries of ${\cal K}_{{\cal N}q'q'_{b}}$.  
  
We now establish two basic existence lemmae for  
(\ref{greq}) and (\ref{greqd}).  
\\
  
{\bf Lemma 1: Fusion of generalized reflection matrices}  
\\  
\\  
 For simplicity we omit the indices $q$ and $q_{b}$.  
  
{\it If} $T$ {\it  is a solution to the reflection equation} (\ref{re})  
{\it then the following objects}:  \be &&{\cal T}^0_{{\cal N}} =  T_{1} R_{21}(\lambda_{1} +\lambda_{2})R_{31}(\lambda_{1}  
+\lambda_{3}) \ldots R_{n1}(\lambda_{1} +\lambda_{n}) T_{2}R_{32}(\lambda_{2} +\lambda_{3})  
\ldots R_{n2}(\lambda_{2} +\lambda_{n}) \non\\  && T_{3} \ldots  T_{k}R_{k+1 k}(\lambda_{k}  
+\lambda_{k+1}) \ldots R_{nk}(\lambda_{k} +\lambda_{n})T_{k+1} \ldots T_{n-1}  
R_{n n-1}(\lambda_{n} +\lambda_{n-1}) T_{n}, \label{sol1} \ee  {\it  are solutions to the set of Generalized Reflection Equations} (\ref{greq}).  
\\  
  
{\bf Lemma 2: Fusion of Dual Generalized reflection matrices}  
\\  
\\  
 Again for simplicity we omit the indices $q'$ and $q'_{b}$.  
  
{\it If} $K$ {\it is a solution to the reflection equation} (\ref{red}) {\it then the following objects}: \be {\cal K}_{{\cal N}}^{0} &=& K_{n}  
M_{n-1} R_{n-1 n}(-\lambda_{n} -\lambda_{n-1}-2\rho)M_{n-1}^{-1}K_{n-1} \ldots  
 \non\\ && K_{k+1}M_{k}R_{kn}(-\lambda_{k} -\lambda_{n}-2\rho) \ldots  
R_{kk+1}(-\lambda_{k} -\lambda_{k+1}-2\rho) M_{k}^{-1} 
K_{k} \non\\ && \ldots K_{2} M_{1} R_{1n}(-\lambda_{1} -\lambda_{n}-2\rho)  \ldots  
R_{12}(-\lambda_{1} -\lambda_{2}-2\rho)M_{1}^{-1} K_{1} \label{sol2} \ee {\it satisfy the dual generalized reflection equation}.  
  
Note that the fusion procedure operates exclusively on auxiliary spaces; the  
quantum spaces $q,q',q_b,q'_{b}$ are untouched.

The proofs are established by a recursion procedure on the total number of  
fused auxiliary spaces in the considered equations:  
 $ n_0 \equiv card ({\cal N} + {\cal M'})$. 
\\  
\\  
{\bf Proof of Lemma 1}  
  
 $\bullet$  The lemma is established by hypothesis for $n_0 = 2$ where  
 necessarily $card ({\cal N}) = card ({\cal M'}) = 1$ and  
(\ref{greq}) reduces to (\ref{re}).  
  
$\bullet$ Assuming now that Lemma 1 is proved up to some value  
 $n_0 \ge  2$ one considers the equations from the set (\ref{greq}) at  
 $n_0 +1$. It is always possible to assume  
that $card ({\cal N})\ge  2$, indeed if not then necessarily  
$card ({\cal M'})\ge  2$  
and by multiplying the equation from  
(\ref{greq}) on the l.h.s.  by $ R_{{\cal N} \bar{\cal  
M'}}(\lambda_{{\cal N}}-\lambda_{\bar{\cal M}'})$ and on the r.h.s. by  
$ R_{\bar{\cal M'} {\cal N}}(\lambda_{{\cal N}}-\lambda_{\bar{\cal M}'})$  
one gets back, after  
exchanging the notations ${\cal N}$ and  ${\cal M'}$, a new form  
of the same equation, this time  with $card ({\cal N}) \ge  2$.  
  
One then particularizes the first index of the ordered set ${\cal N}$ as $1$  
and rewrites the equation as (denoting the set ${\cal N} - \{ 1 \}$ as  
${\cal N}^{-}$):  
\be  
&& {\cal T}_{1}(\lambda_{1})\  
 R_{{\cal N}^{-}{1}}(\lambda_{ 1}+\lambda_{ {\cal N}^{-}})\  
 {\cal T}_{{\cal N}^{-}}(\lambda_{{\cal N}^{-}})\  
 R_{{\cal M}'1}(\lambda_{1}+\lambda_{ {\cal M}'})\  
 R_{{\cal M}'{\cal N}^{-}}(\lambda_{ {\cal  
N'}}+\lambda_{ {\cal M}'})\ \non\\ &&{\cal T}_{{\cal M}'}(\lambda_{{\cal  
M}'})\ R_{\bar {\cal N}^{-} {\cal M}'}(-\lambda_{\bar  {\cal  
N}^{-}}+\lambda_{{\cal M}'})\  
 R_{1 {\cal M}'}(-\lambda_{1}+\lambda_{{\cal M}'})  
 \non\\  &=&  
 R_{ {\cal M}'\bar{\cal  
N}^{-}}(-\lambda_{\bar  {\cal N}^{-}}+\lambda_{ {\cal M}'})\  
R_{ {\cal M}' 1}(-\lambda_1+\lambda_{ {\cal M}'})\  
 {\cal T}_{{\cal M}'}(\lambda_{{\cal M}'})\  
R_{1  {\cal  
M}'}(\lambda_1+\lambda_{ {\cal M}'}) \non\\ &&  
 R_{{\cal N}^{-}  {\cal  
M}'}(\lambda_{ {\cal N}^{-}}+\lambda_{ {\cal M}'})\  
{\cal T}_{1}(\lambda_{1})\  
 R_{{\cal N}^{-}{1}}(\lambda_{ 1}+\lambda_{ {\cal N}^{-}})\  
 {\cal T}_{{\cal N}^{-}}(\lambda_{{\cal N}^{-}}). \label{greq3} \ee  One now
 establishes validity of this equality by successive  
operations on the l.h.s. of (\ref{greq3}):  
  
{\bf 1} using the (already proved by recursion hypothesis) exchange  
relation for index sets ${\cal N}^{-}$ and  ${\cal M'}$.  
  
{\bf 2} using a fused Yang--Baxter equation: \be R^{+}_{{\cal  
N}^{-}1}\ R^{+}_{{\cal M'}1}\  R^{-}_{{\cal M'}\bar{\cal N}^{-}} =  
R^{-}_{{\cal M'} \bar{\cal N}^{-}}\  R^{+}_{{\cal M'}1}\  R^{+}_{{\cal  
N}^{-}1} \label{fyb1} \ee where the compact notations $R^{\pm}$  
are self-explanatory.  
  
{\bf 3} using a second fused Yang-Baxter equation:  
 \be R^{+}_{{\cal N}^{-}1}\ R^{+}_{{\cal N}^{-}{\cal M'}}\ R^{-}_{1 {\cal M'}}  
= R^{-}_{1 {\cal M'}}\ R^{+}_{{\cal N}^{-}{\cal M'}}\ R^{+}_{{\cal N}^{-}1}  
\label{fyb2} \ee  
  
{\bf 4} using the (already proved by recursion hypothesis)  
exchange relations for sets  $\{ 1 \}$ and  ${\cal M'}$.  
  
Note that the fused Yang-Baxter equations (\ref{fyb1}), (\ref{fyb2}) are immediate  
consequences of the coproduct structure of the ordinary Yang-Baxter equation.  
  
This establishes Lemma 1.  
  
{\bf Proof of  Lemma 2}  
  
The proof is similar, successively applying to the r.h.s.  
of (\ref{greqd}) decoupled as in (\ref{greq3})  
  
{\bf 1} the reflection equation for sets ${\cal N}^{-}$ and  ${\cal M'}$.  
  
{\bf 2} the dual fused Yang Baxter equation \be R^{+}_{{\cal  
N}^{-} \bar{\cal M'}}\ R^{--}_{1 \bar{\cal M'}}\ R^{--}_{1 \bar {\cal  
N}^{-}} = R^{--}_{1 \bar {\cal N}^{-}}\ R^{--}_{1 \bar{\cal M'}}\  
 R^{+}_{{\cal N}^{-} \bar{\cal M'}}  
\label{dfyb1}  
\ee where the notation $R^{--}$ is also self-explanatory from (\ref{greqd})  
  
{\bf 3} the dual fused Yang Baxter equation  
\be R^{-}_{\bar{\cal M'}1}\  R^{--}_{ \bar{\cal M'}\bar{\cal N}^{-}}\ R^{--}_{1 \bar {\cal N}^{-}}  
=  R^{--}_{1 \bar {\cal N}^{-}}\ R^{--}_{ \bar{\cal M'}\bar{\cal N}^{-}}\  
 R^{-}_{\bar{\cal M'}1}  
\label{dfyb2} \ee  
  
{\bf 4} the reflection equation for $\{1 \}$ and  ${\cal M'}$,  
  
{\bf 5} and at several places the commutation relation (\ref{crosscom2}). 
 
Both lemmas are thus established, hence the set of solutions of  
(\ref{greq}), (\ref{greqd}) is not empty.  
  
Notice that the fused solution  ${\cal K}^0_{{\cal N}}$ has a  
structure similar to  ${\cal T}^0_{{\cal N}}$, but with a reversed order  
of the auxiliary spaces; in particular one  has an identification of formal solutions  
as:  
${\cal K}^0_{{\cal N}}(\lambda_{{\cal N}}) \equiv  
 {\cal T}^0_{{\cal N}}(-\lambda_{{\cal N}}-\rho)^{t_{{\cal N}}}$.  
 This is  expected because of the form of the equations (\ref{greq})  
 and (\ref{greqd}): equation (\ref{greq}) is formally  
 the ``transposed'' of equation (\ref{greqd}).  The treatment
of the general case will follows from the fact that fused YB
equations (\ref{fyb1}), (\ref{fyb2}), (\ref{dfyb1}), (\ref{dfyb2}) 
take identical forms in the general case once one replaces \be
 R_{12}(\lambda_1 - \lambda_2) \rightarrow  A, ~~ R_{21}(\lambda_1 - \lambda_2)
 \rightarrow D, ~~ R^{+} \rightarrow C = B^{\pi}. \ee This replacement will provide us with a ``dictionary'' between the reflection
algebra structures and the general quadratic structures.
\\  
\\  
We now establish a very important ``dressing'' or invariance  
property of general solutions to the fused equations. It is described by:

{\bf Proposition 1: dressing of solutions}  
  
{\bf 1a}
{\it Given a set of solutions to the dual generalized reflection equations}  
(\ref{greqd}) {\it one obtains a new set of solutions by multiplying}  
${\cal K}_{{\cal N}, q}$ {\it on the left with operators of the form}  
${\cal Q}_{\cal N}$ {\it acting on the sole auxiliary spaces, and such that  
 for all sets} ${\cal M'}$ {\it disjoint from} ${\cal N}$ {\it one has}:  
\be  
[ {\cal Q}_{\cal N},\ R_{{\cal M'}{\cal N}}(\lambda_{{\cal N}} \pm \lambda_{{\cal M}'} )]= 0,  
 ~~ [ {\cal Q}_{\cal N},\ R_{\bar {\cal N}{\cal M}'}(-\lambda_{{\cal N}} \pm \lambda_{{\cal M}'} )] = 0.  
\label{dress}  
\ee
{\it A similar property will hold for} (\ref{greq}){\it under right multiplication by}  ${\cal Q}_{{\cal N}}$.

{\bf 1b}
{\it Given a set of solutions to the dual generalized reflection equations}  
(\ref{greqd}) {\it one obtains a new set of solutions by multiplying}  
${\cal K}_{{\cal N}, q}$ {\it on the right with operators of the form}  
${\cal S}_{\cal N}$ {\it acting on the sole auxiliary spaces, and such that  
 for all sets} ${\cal M'}$ {\it disjoint from} ${\cal N}$ {\it one has}:  
\be  
[ {\cal S}_{\cal N},\ R_{{\cal N}{\cal M'}}( \lambda_{{\cal N}} \pm \lambda_{{\cal M}'} )]= 0,  
 ~~ [ {\cal S}_{\cal N},\ R_{{\cal M}'\bar {\cal N}}(- \lambda_{{\cal N}} \pm \lambda_{{\cal M}'} )] = 0.  
\label{dress1}  
\ee
{\it A similar property will hold for} (\ref{greq})
{\it under left multiplication by}  ${\cal S}_{{\cal N}}$.

The two conditions in (\ref{dress} and in (\ref{dress1})
are equivalent by unitarity of fused $R$-matrices.

The proof is immediate from the form of  (\ref{greq}),(\ref{greqd}) and the use
of (\ref{crosscom2}).

We shall discuss for simplicity
from now on only the case of dressing operators of type 
${\cal Q}_{{\cal N }}$. Notice that in particular the 
following product of $\check R$'s defined by (\ref{dress0}): \be  
{\cal Q}_{{\cal N}}=\check R_{12}(\lambda_{1} -\lambda_{2}) \ldots \check R_{n-1 n}(\lambda_{n-1} -\lambda_{n}) \label{sol3} \ee realizes such a dressing. As indicated in section 2.1. the construction of $\check R$ as $P_{12}R_{12}$ is problematic when considering spectral--parameter dependent $R$ matrices where the permutation operator must act on a loop space $(V_{1} \otimes C(\lambda_{1})) \otimes (V_{2} \otimes C(\lambda_{2}))$. One may however establish at the level of formal series a
 better defined solution as: 
  
{\bf Proposition 2}
\be \check R = \delta (\lambda_{1},\ \lambda_{2})\Pi_{12} R_{12} \label{check} \ee {\it realizes}  (\ref{dress0}) {\it as formal series, where $\Pi_{12}$ is the operator exchanging vector spaces $V_{1}$ ,$V_{2}$}.

Here the $\delta$ distribution is of course defined as a formal series
\be \delta(\lambda_{1},\ \lambda_{2}) =\sum_{n \in {\bf Z}}(\lambda_{1} / \lambda_{2})^{n}. \ee
The proof follows from the formal series identity:
\be
\delta (\lambda_{1},\ \lambda_{2})\ f(\lambda_{1})\ g(\lambda_{2}) = f(\lambda_{2})\ g(\lambda_{1})\ \delta (\lambda_{1},\ \lambda_{2}) \label{d} \ee
for any functions $f$, $g$ with an assumed Laurent series expansion defined as  \be f(\lambda) = \sum_{-\infty<m_{0}<m} f_{m} \lambda^{m}, ~~ g(\lambda) = \sum_{-\infty < p_{0} < p} g_{p} \lambda^{p}. \ee The explicit proof follows from comparing both sides of (\ref{d}) as formal series: the coefficients of $\lambda_{1}^{a}\lambda_{2}^{b}$ on both sides are \be\sum_{m,p|m+p=a+b} (f_{m} g_{p})\ee which is a finite sum by the Laurent series hypothesis on $f$, $g$ ($-\infty <m_{0}$, $-\infty < p_{0}$). Provided that one deals with $R$ matrices expandable as Laurent series, one then immediately deduces that (\ref{check}) satisfies (\ref{dress0}), and hence (\ref{sol3}) satisfies (\ref{dress}).

As in the Yang--Baxter type algebras, having determined the proper
 generalized exchange relations, we are now in a position  to build  
commuting traces. We establish the fundamental  
\\  
\\  
{\bf Theorem}  
  
{\it Let} ${\cal K}_{{\cal N}}$ {\it  be a set of solutions to the  
dual generalized reflection equations} (\ref{greqd}), {\it  acting  
on the auxiliary spaces labelled by} ${\cal N}$, {\it the quantum  
space labelled by} $q'$, {\it and a possible boundary space  
labelled by} $q_{b}'$;  
  
{\it Let} ${\cal T}_{{\cal N}}$ {\it be a set of solutions of the generalized  
 reflection equations} (\ref{greq}) {\it acting on the tensor product of  
 the auxiliary spaces labelled by} ${\cal N}$, {\it the quantum space labelled by}  $q$,  
 {\it and the boundary space labelled by} $q_{b}$.  
  
{\it The following trace operators acting on the quantum space}
 $q \otimes q_b \otimes q' \otimes q'_b$:  
\be  
H_{{\cal N}} = Tr_{{\cal N}} \Big ( {\cal K}_{{\cal N}}^{*}(\lambda_{{\cal N}}){\cal T}_{{\cal N}}(\lambda_{{\cal N}}) \Big ),  
\label{tracer} \ee  
{\it where} ${\cal K}_{{\cal N}}^{*}= {\cal K}_{{\cal N}}^{t_{q'q'_{b}}}$,  
{\it  build  a family of mutually commuting operators}:  
\be \Big [H_{{\cal N}}, \ H_{{\cal M}'} \Big ] = 0. \label{comrel} \ee  
{\bf Proof}  
  
The proof is identical to the proof of (\ref{com0})
 (see \cite{sklyanin})  
thanks to the fusion  relations of Lemma 1 and 2; it is however
 worth being given in  
detail. Quantum indices will again be dropped for simplicity.  
  
One starts from the product $H_{{\cal M}'}H_{{\cal N}}$;
 the idea is to pass the unprimed  
indices of the product through the primed indices. For
 this purpose we first act by the  
partial transposition $t_{{\cal M}'}$, yielding:
\be H_{{\cal M}'}H_{{\cal N}}  
 &=& \tr_{{\cal M}'} {\cal  K}_{{\cal M}'}^{*}(\lambda_{{\cal M}'})\ {\cal T}_{{\cal M}'}  
 (\lambda_{{\cal M}'}) \tr_{{\cal N}} {\cal K}_{{\cal N}}^{*}(\lambda_{{\cal N}})\  
 {\cal T}_{{\cal N}}(\lambda_{{\cal N}}) \non \\ &=& \tr_{{\cal M}'}  
 {\cal K}_{{\cal M}'}^{*}(\lambda_{{\cal M}'})^{t_{{\cal M}'}}\ {\cal T}_{{\cal M}'}
 (\lambda_{{\cal M}'})^{t_{{\cal M}'}}\tr_{{\cal N}}  {\cal K}_{{\cal N}}^{*}  
 (\lambda_{{\cal N}})\ {\cal T}_{{\cal N}}(\lambda_{{\cal N}}) \non \\  
 &=& \tr_{{\cal M}'{\cal N}}  {\cal K}_{{\cal M}'}^{*}(\lambda_{{\cal M}'})^{t_{{\cal M}'}}\  
   {\cal K}_{{\cal N}}^{*}(\lambda_{{\cal N}})\ {\cal T}_{{\cal M}'}(\lambda_{{\cal M}'})^{t_{{\cal M}'}}\  {\cal T}_{{\cal N}}(\lambda_{{\cal N}})\,.  
\ee  
We now use the crossing-unitarity of the $R$ matrix  
\be M_{{\cal M}'}^{-1}\ R_{\bar {\cal M}'\bar {\cal N}}(-\lambda_{\bar {\cal M}'}-\lambda_{\bar {\cal N}} -2 \rho)^{t_{{\cal M}'}}\ M_{{\cal M}'}\ R_{{\cal N}{\cal M}'}(\lambda_{{\cal M}'}+\lambda_{{\cal N}})^{t_{{\cal M}'}} = Z(\lambda_{{\cal M}'}+\lambda_{{\cal N}})\,, \label{1} \ee  
where $Z$ is just a function of $\lambda$'s. Then the product $H_{{\cal M}'}H_{{\cal N}}$ becomes \be  
&&Z^{-1}(\lambda_{{\cal M}'}+\lambda_{{\cal N}}) \tr_{{\cal M}'{\cal N}}  
 {\cal K}_{{\cal M}'}^{*}(\lambda_{{\cal M}'})^{ t_{{\cal M}'}}\  
 {\cal K}_{{\cal N}}^{*}(\lambda_{{\cal N}})\  
M_{{\cal M}'}^{-1}\  
R_{{\cal N} {\cal M}'}(-\lambda_{ {\cal M}'}-\lambda_{ {\cal N}} -2 \rho)^{t_{{\cal N}}}\ M_{{\cal M}'}  
\non\\ && \times  R_{{\cal N} {\cal M}'}(\lambda_{{\cal M}'}+\lambda_{{\cal N}})^{t_{{\cal M}'}}\  {\cal T}_{{\cal M}'}(\lambda_{{\cal M}'})^{t_{{\cal M}'}}\  {\cal T}_{{\cal N}}(\lambda_{{\cal N}}) \non \\&&  
=  Z^{-1}(\lambda_{{\cal M}'}+\lambda_{{\cal N}}) \tr_{{\cal M}'{\cal N}}  
\Big ({\cal K}_{{\cal M}'}^{*}(\lambda_{{\cal M}'})^{t_{{\cal M}'}}\ M_{{\cal M}'}^{-1}\ R_{{\cal N}{\cal M}'}(-\lambda_{{\cal M}'}- \lambda_{{\cal N}} -2 \rho)\ M_{{\cal M}'}\  {\cal K}_{{\cal N}}^{*}(\lambda_{{\cal N}})\Big )^{t_{{\cal N}}}  
 \non\\ && \times \Big ({\cal T}_{{\cal M}'}(\lambda_{{\cal M}'})\  R_{{\cal N}{\cal M}'}(\lambda_{{\cal M}'}+\lambda_{{\cal N}})\ {\cal T}_{{\cal N}}(\lambda_{{\cal N}})\Big  )^{t_{{\cal M}'}} \non\\&& = Z^{-1}(\lambda_{{\cal M}'}+\lambda_{{\cal N}}) \tr_{{\cal M}'{\cal N}}  
\Big ( K_{{\cal M}'}^{*}(\lambda_{{\cal M}'})^{t_{{\cal M}'}}\ M_{{\cal M}'}^{-1}\ R_{ {\cal N} {\cal M}'}(-\lambda_{{\cal M}'} -\lambda_{{\cal N}} -2 \rho)\  M_{{\cal M}'}\  
 K_{{\cal N}}^{*}(\lambda_{{\cal N}}})^{t_{{\cal N}}}\Big )^{t_{{\cal M}'{\cal N}} \non\\ && \times {\cal T}_{{\cal M}'}(\lambda_{{\cal M}'})\  R_{{\cal N} {\cal M}'}(\lambda_{{\cal M}'}+\lambda_{{\cal N}}) {\cal T}_{{\cal N}}(\lambda_{{\cal N}})\,.  
\ee  
Using the unitarity of the $R$ matrix \be R_{\bar {\cal N}\bar {\cal M}'}(-\lambda_{\bar {\cal M}'}+\lambda_{\bar {\cal N}})\ R_{{\cal M}'{\cal N}}(\lambda_{{\cal M}'}-\lambda_{{\cal N}}) =  
 Z(\lambda_{{\cal M}'}-\lambda_{{\cal N}})\,,\label{2}\ee  
we obtain the following expression for the product  
\be && Z^{-1}(\lambda_{{\cal M}'}-\lambda_{{\cal N}})Z^{-1}(\lambda_{{\cal M}'}+\lambda_{{\cal N}}) \tr_{{\cal M}'{\cal N}} \Big ( {\cal K}_{{\cal M}'}^{*}(\lambda_{{\cal M}'})^{t_{{\cal M}'}}\ M_{{\cal M}'}^{-1}\ \non\\ && \times R_{{\cal N} {\cal M}'}(-\lambda_{{\cal M}'} -\lambda_{ {\cal N}} -2 \rho)\ M_{{\cal M}'}\  {\cal K}_{{\cal N}}^{*}(\lambda_{{\cal N}})^{t_{{\cal N}}}\Big )^{t_{{\cal M}'{\cal N}}}\non \\&& \times R_{{\cal N}\bar {\cal M}'}(-\lambda_{\bar {\cal M}'}+\lambda_{{\cal N}})\ R_{{\cal M}'\bar {\cal N}}(\lambda_{{\cal M}'}-\lambda_{\bar {\cal N}}) {\cal T}_{{\cal M}'}(\lambda_{{\cal M}'})\  
  R_{{\cal N} {\cal M}'}(\lambda_{{\cal M}'}+\lambda_{{\cal N}})\  
{\cal T}_{{\cal N}}(\lambda_{{\cal N}})  
\non \\ &&=  Z^{-1}(\lambda_{{\cal M}'}-\lambda_{{\cal N}}) Z^{-1}(\lambda_{{\cal M}'}+\lambda_{{\cal N}}) \tr_{{\cal M}'{\cal N}} \Big  (R_{{\cal M}'\bar {\cal N}}(-\lambda_{{\cal M}'} + \lambda_{\bar {\cal N}})\   {\cal K}_{{\cal M}'}^{*}(\lambda_{{\cal M}'})^{t_{{\cal M}'}}\ \non\\ && \times M_{{\cal M}'}^{-1}\  R_{{\cal N}  {\cal M}'}(-\lambda_{ {\cal M}'}-\lambda_{{\cal N}} -2 \rho)\ M_{{\cal M}'}\   {\cal K}_{{\cal N}}^{*}(\lambda_{{\cal N}})^{t_{{\cal N}}}\Big  )^{t_{{\cal M}'{\cal N}}}\non \\&& \times R_{{\cal M}'\bar {\cal N}}(\lambda_{{\cal M}'}- \lambda_{\bar {\cal N}})\ {\cal T}_{{\cal M}'}(\lambda_{{\cal M}'})\  R_{{\cal N} {\cal M}'}(\lambda_{{\cal M}'} + \lambda_{{\cal N}} )\ {\cal T}_{{\cal N}}(\lambda_{{\cal N}})\,.\ee  
One now recognizes in the r.h.s.: {\bf 1}. the r.h.s. of the exchange equation (\ref{greq}); {\bf 2}. the full transposition under $t_{{\cal N}{\cal M}'}t_{q'q'_{b}}$ of the l.h.s. of the exchange relation (\ref{greqd}) (recall (\ref{p1})). The specific forms of unitarity (\ref{2}) and crossing--unitarity (\ref{1}) are crucial in yielding this form. This hints at a close connection between the crossing--unitarity properties and the duality ``transformation'' between (\ref{greq}) and (\ref{greqd}).  
  
Now, with the help of equations (\ref{greq}), (\ref{greqd}), (\ref{1}) and (\ref{2}), and  
by repeating the previous steps in a reverse order we establish  
that the last expression is indeed $H_{{\cal N}}H_{{\cal M}'}$. This  
concludes the proof of the commutativity relation (\ref{comrel}).  
  
Note that it is needed at the very beginning of  this proof to use mutual  
commutation of the matrix elements of ${\cal K}$ and ${\cal T}$. We then
emphasize  
that the above proof is  
valid as long as neither the left and right boundary spaces $q_b$ and
$q'_b$, nor the quantum spaces $q$ and $q'$ ``talk'' to each other.
   
\subsection{Dressed quantum traces}

The dressing procedure of fused matrices was defined by {\bf Proposition 1}. It generically leads to traces of the form $Tr({\cal Q}{\cal K}^{0}{\cal S}{\cal T}^{0})$.
We shall again consider only at this stage dressing by objects of type 
${\cal Q}$ leading to traces of the form $Tr({\cal Q}{\cal K}^0{\cal T})$
so as to avoid complicating the discussion. 

This construction is needed to get non-trivial commuting traces.  
Indeed, it is easy to prove:  
  
{\bf Proposition 3: Factorized quantum traces}  
  
{\it Operators built from the basic fused solutions}  
 (\ref{sol1}), (\ref{sol2}) {\it decouple as}  
$H_{{\cal N}} = Tr_{{\cal N}} 
\Big ({\cal K}^{0 *}_{{\cal N}}{\cal T}^0_{{\cal N}}\Big ) 
= \Big (Tr_{1}K_{1}^{*} T_1\Big )^{card {\cal N}}$.  
  
The proof is achieved by successive use of partial transpositions with  
respect to successive indices of ${\cal N}$ from $1$ to $n$ as defined in  
 (\ref{sol1}), (\ref{sol2}). These partial transpositions systematically  
bring together fused transposed $R$ matrices which cancel each other by  
fused crossing relations (\ref{cross2}). Hence the complete elimination of the  
$R$ matrices, yielding
 $ Tr_{{\cal N}} \Big ({\cal K}_{{\cal N}}^{0*}{\cal T}^0_{{\cal N}}\Big )= Tr \Big ((K_{n} \ldots K_{1})^{*} T_{1} \ldots T_{n}\Big )$,  then by using the antimorphism property of the $*$ operation one gets
 $Tr \Big (K_{1}^{*}
 \ldots K_{n} ^{*} T_{1} \ldots T_{n}\Big )=\Big (Tr K_{1}^{*} T_{1}\Big )^{n}$. 
Such was also the more obvious case for usual Yang Baxter algebras  
(\ref{trace}), (\ref{sol}).  
  
Therefore additional non-trivial dressing operators are required  
to get non trivial traces.  As indicated by {\bf Proposition 2} , examples of  
such objects are already formally available as products of neighboring-indexed 
$\delta(\lambda_{a},\ \lambda_{a+1})\ \Pi_{a,a+1}\ R_{a,a+1}$. This particular choice of dressing deserves a more detailed discussion. It provides actually an interesting classical limit. We here assume that the classical limits are defined in the usual way, i.e.
\be R_{12} (\lambda_1 - \lambda_2) = {\bf 1} \otimes {\bf 1} + \hbar r_{12}(\lambda_1 - \lambda_2); ~~T(\lambda)
= t(\lambda)+ o(\hbar);~~K(\lambda)
= {\bf 1} + o(\hbar) \label{classlim} \ee
One then easily shows that:

{\bf Proposition 4}

{\it The coefficients of the operators} $ H_{{\cal N}} = Tr_{{\cal N}}
 \Big ( {\cal K}_{{\cal N}}^{*}(\lambda_{{\cal N}}){\cal T}_{{\cal N}}(\lambda_{{\cal N}}) \Big)$ {\it expanded as a multiple formal series with general term}
$\lambda_1^{a_1} \cdots \lambda_n^{a_n}, n = ${\it card}$({\cal N})$ 
{\it are identical, when} $\hbar $
{\it goes to zero, to the coefficients of the classical trace}    
$ h_{n} = Tr(t^n(\lambda))$ 
{\it expanded as a formal series with general term } $\lambda^m$
{\it with the identification} $ m = \sum_{i=1}^{n} a_i$.

{\bf Proof} 

The properties of the permutation operators
$\Pi_{a,a+1}$  inside the multiple trace reduce it to a single trace
of direct products, once both $R$ and $K$ factors are taken to identity. 
The formal series expansion of the product of $\delta$ distributions then yields
the suggested identification at $m = \sum_{i=1}^{n} a_i$
of the coefficients of the two formal series. Again we assume that $t(\lambda)$ can be formally expanded as Laurent series.

Since the classical limit of these particular dressed traces yields the 
classical trace of a power of the classical Lax matrix we are justified, as
was the case in \cite{MAI1}, in denoting these operators as ``quantum traces'',
a priori {\it algebraically distinct} from operators obtained from powers of 
the commuting traces of the quantum Lax matrix. 
 
Let us conclude this brief discussion by a general remark.
The problem of  
constructing generators such that $\Big [{\cal Q}_{\cal N},\ R_{\cal NM}\Big ]
 =0$ for  
general fused $R$ matrices is in any case related to the  
understanding of the underlying coproduct structure and universal  
algebra. In this context ${\cal Q}$ objects may be constructed generically  
as coproducts of central elements \footnote{We are indebted to Daniel Arnaudon  
for this suggestion}.

\section{Conclusion and prospective}  
  
In order to define lines of future investigation it is important to  
characterize the basic steps of the procedure described in the previous  
section. It crucially depends on four fundamental features:  
  
{\bf Step 1 -} Existence of an algebraic reflection-like structure (\ref{re})  
with the notion of a dual structure (\ref{red}); associated unitarity and crossing  
relations, and dual trace formula generating commuting objects as in (\ref{trace1}).
   
{\bf Step 2 -} Existence of mutually  
consistent fusion procedures for both algebraic structures as described in  
Lemma 1 and 2.  
  
{\bf Step 3 -} Dressing of  fused solutions by commuting (fused) operators  
${\cal Q}$ on the auxiliary spaces, such as characterized by Proposition 1.  
  
{\bf Step 4 -} One then combines 1, 2 and 3 to get commuting traces  
by products of fused solutions ${\cal T}$ and ${\cal K}$.  
  
This now indicates several directions of investigation:

{\bf Extension of the procedure to general reflection algebras of 
(}\ref{braid}{\bf) type}  
 
Interest in this generalization stems in particular from the  
occurence in physical systems of relevant examples of braided YB  
algebras: the structure matrices $A, B, C, D$ are in general not  
independent or free. For instance the reflection algebra itself is  
characterized by $A_{12} =D_{21} = R_{12}(\lambda_{1}  
-\lambda_{2})$ and $B_{12} =C_{21} = R_{12}(\lambda_{1}  
+\lambda_{2})$. Different choices of constraints realize different  
algebraic objects. Such are for instance the  
reflection-transmission algebras \cite{MSR} where  $A, B, C, D$  
are given by one single quantum $R$  
 matrix depending however on two independent spectral parameters instead of  
one single combination $\lambda_1 \pm \lambda_2$;  the non-ultralocal  
monodromy matrix algebra for mKdV equations \cite{Ku,Ro} where the $B$  
matrices are diagonal c-number  matrices; and  cases where the $R$ matrix  
does not necessarily satisfy the crossing symmetry (\ref{cross}), but a  
more general form of crossing (see e.g. \cite{doikou2}).  
 
The extension will be undertaken \cite{ADN} by following Steps 1 to 4  
as previously defined.

 Step 1 was actually realized in \cite{MAI2} where a dual  
reflection algebra was defined.  
It is indeed easy to check that it gives back (\ref{red}) in our case  
once one uses unitarity, crossing and transposition relations  
(\ref{cross}) and (\ref{transp}). The commuting traces are again of the form  
(\ref{trace1}). However in \cite{MAI2} the matrices ${\cal K}$ are taken to be  
pure c-number objects, not quantum operators; hence it remains to prove consistency of  
the definition for quantum ${\cal K}$ matrices.

The next problem is to realize Step 2, defining fusion 
procedures for quadratic algebras generalizing  Lemma 1 and 2. 
Fusion procedures are obtained easily  
once a coproduct structure has been formulated; however consistent
fusion procedures may exist at a represented level without
being direct consequences of underlying universal coproduct structures:
for instance it seems to be the case for the formulation given
here\footnote{as follows from our enlightening discussions with P.P. Kulish}.

By contrast another fusion procedure exists, compatible 
with (\ref{re}), defined in \cite{KDM} as following from a universal coproduct
structure. Connections between these two fusion
processes will be discussed in detail in the general context. Briefly, the
basic ingredient in this connection is a ``braiding'' matrix 
${\cal F}_{{\cal N}}$, i.e. acting on fused R-matrices as:
\be {\cal F}_{{\cal N}} R_{ {\cal N} {\cal M}'} = R_{{\bar {\cal N}} {\cal M}'}{\cal F}_{{\cal N}}  \label{braid2} \ee This procedure can naturally also be used to define quantum traces from (\ref{tracer}). The comparison between the two constructions
will be left for the general discussion. 
      
It is thus possible to
define {\it several}  fusion procedures and trace formulae,
depending on the considered subclasses of quadratic algebras, characterized  
in particular by different restrictions on the structure matrices  
$ A, B, C, D$.

In view of physical applications such
as already mentioned, it is in fact very interesting to  
examine which restrictions on $A$, $B$, $C$, $D$ are compatible with the duality  
construction and how they relate to particular crossing and unitarity conditions.  
This in turn will need a careful analysis  
of the discrete (possibly infinite) symmetry groups of the braided YB  
equations analogous to the studies conducted e.g. in \cite{BMV}.

Finally we have already discussed the extensions of Step 3.

{\bf  Extensions to quadratic algebras of dynamical type}  
  
Here the matrices $R$ and $K$  should depend on some extra  
``dynamical'' parameters identified as coordinates on the dual of  
some (possibly non-abelian, see \cite{Ping}) Lie algebra, and the  
exchange relations take a form along the lines of the dynamical  
Yang--Baxter equation described in \cite{GFN,Ping}. The first  
problem at this time is to actually construct such extensions.  
Some examples are currently being considered  
\cite{Z, ER}. Quantum traces formulas may then be obtained by  
suitably manipulating our construction along the lines in  
\cite{ABB} for dynamical Yang--Baxter equation, so as to  
incorporate the shifted dynamical parameter. Again one essentially  
needs to define consistent fusion procedures for such algebras,  
and consistent unitarity and crossing properties.  
  
In particular, if one obtains solutions of the example of dynamical quantum reflection equation defined in \cite{Z}, which incorporates as a structure matrix Felder's dynamical $R$ matrix \cite{GFN}, one expects to be able to  
define generalizations of the Ruijsenaars--Schneider Hamiltonians \cite{RS}  
obtained in \cite{ABB}. Such generalizations, incorporating non--trivial  
 reflection  
algebras, may be connected with the quantum Hamiltonians presented  
e.g. in \cite{Koo, VD, HK, AR} as the construction \cite{HK} suggests,  
however more general integrable Hamiltonians may also arise.  
  
{\bf Acknowledgments}  
  
We are grateful to D. Arnaudon, L. Frappat and E. Ragoucy for  
helpful discussions. We are also indebted to J.M. Maillet and P.P. Kulish for  
valuable comments and discussions. A.D. is supported by the TMR  
Network ``EUCLID''; ``Integrable models and applications: from  
strings to condensed matter'', contract number  
HPRN--CT--2002--00325. J.A. thanks as ever LAPTH Annecy for kind support.

\end{document}